\documentclass[12pt]{amsart}
\usepackage{amssymb}
\newtheorem{prop}{Proposition}
\newtheorem{lemma}{Lemma}
\newtheorem{theorem}{Theorem}
\newtheorem{corollary}{Corollary}
\newtheorem{remark}{Remark}
\newtheorem{definition}{Definition}

\newcommand{\RR}{{\mathbb R}}

\newcommand{\ZZ}{{\mathbb Z}}

\begin{document}

\author{V.A.~Vassiliev}

\address{Steklov Mathematical Institute of Russian Academy of Sciences} \email{vva@mi-ras.ru}

\title[Average intersection number of plane curves]{Average intersection number of trigonometric plane curves \\ in $L_2$ and $W_2^r$ statistics }

\begin{abstract}
The average value of the number of intersection points of two plane curves given by trigonometrical polynomial maps of degree $N$ and bounded $L_2$- or $W_2^r$-  norms is calculated
\end{abstract}



\maketitle

\section{Main result}

Given a natural number $N$, consider the space of maps $S^1 \to \RR^2$ given by a pair of Fourier polynomials of degree $\leq N$,
\begin{equation} 
\label{on} \left\{
\begin{array}{ccl}
x(\varphi) & =  a_0 & +a_1 \cos \varphi + \dots + a_N \cos N \varphi + \\ & & + 
b_1 \sin \varphi + \dots + b_N \sin N \varphi  \ , \\
y(\varphi) & = \tilde a_0 & +\tilde a_1 \cos \varphi + \dots + \tilde a_N \cos N \varphi + \\
& & + \tilde b_1 \sin \varphi + \dots + \tilde 
b_N \sin N \varphi \ .
\end{array}
\right.
\end{equation}
Supply this space with the usual $L_2$ Euclidean structure:
the dot product of two curves $(x_1(\varphi),y_1(\varphi))$ and $(x_2(\varphi),y_2(\varphi))$ is equal to
$$\frac{1}{\pi}\int_{0}^{2\pi} \left(x_1(\varphi)x_2(\varphi) + y_1(\varphi)y_2(\varphi)\right) d \varphi.$$ 
In other words, if another curve is given by\begin{equation}
\label{twoo}
\left\{
\begin{array}{ccl}
x(\varphi) & =c_0 & +c_1 \cos \varphi + \dots + c_N \cos N \varphi + \\
& & + d_1 \sin \varphi + \dots + d_N \sin N \varphi  \ , \\
y(\varphi) & = \tilde c_0 & +\tilde c_1 \cos \varphi + \dots + \tilde c_N \cos N \varphi + \\
& & + \tilde d_1 \sin \varphi + \dots + \tilde 
d_N \sin N \varphi \ , 
\end{array}
\right.
\end{equation}
then its dot product with (\ref{on}) is equal to
$$2\left(a_0c_0 +\tilde a_0\tilde c_0 \right)+ \sum_{j=1}^N\left(a_jc_j+b_jd_j+\tilde a_j\tilde c_j + \tilde b_j\tilde d_j\right).$$
 
The definitions of this space and its Euclidean structure are invariant under the rotations of the plane $\RR^2$.

Denote by $B(N,R)$ the ball of radius $R$ centered at the origin in the $2(2N+1)$-dimensional Euclidean space thus obtained. 
Our main space is $B^2(N,R)$, the space of pairs of such maps. It is also a Euclidean manifold, with the usual Euclidean structure of the direct sum of two Euclidean spaces. For almost all (i.e. constituting a subset of full measure) points of $B^2(N,R)$ the number of the intersection points of the corresponding two plane curves is unambiguous. Denote by $\tilde A_i(N,R)$ the volume of the set of such pairs with exactly $i$ intersection points, and by $A_i(N)$ the corresponding relative volume, that is, $\tilde A_i(N,R)$ divided by the volume of entire $B^2(N,R)$. Obviously, $\tilde A_i(N,R) = A_i(N) =0$ for all odd $i$. 
From the dilation of the plane $\RR^2$, it is easy to see that the numbers $A_i(N)$ indeed do not depend on $R$, so in what follows we assume $R=1$. The volume of the corresponding space $B^2(N,1)$ is equal to
\begin{equation}
\label{totvol}
\left(\frac{\pi^{2N+1}}{(2N+1)!} \right)^2 .
\end{equation}

\begin{theorem}
\label{mth}
The average number 
\begin{equation}
\label{meean}
\sum_{i=0}^\infty i \cdot A_i(N)
\end{equation}
of intersection points of pairs of curves from $B(N,R)$ is equal to 
\begin{equation}
\label{mfor}
\frac{2^{8N+3}((2N)!)^4}{((4N+1)!)^2} \left(1+4 + \cdots + N^2\right).
\end{equation}
\end{theorem}  

\begin{corollary}
\label{maincor}
The number $($\ref{meean}$)$ grows as $\frac{\pi}{3}N^2+ O_{N \to \infty}(N).$ \hfill $\Box$
\end{corollary}

\begin{corollary}
\label{corint}
The same formula $($\ref{mfor}$)$ describes the average number of intersection points of arbitrary two curves of the form $($\ref{on}$)$ with respect to any integrable probability distribution on the space of pairs of such curves, which depends only on the maximum of $L_2$-norms of two curves. \hfill $\Box$
\end{corollary}

A proof of theorem \ref{mth} occupies the sections \ref{fred}---\ref{calin} below. \medskip

In \S \ref{sob} we give an analogous calculation for curves with bounded Sobolev $W_2^r$ norms, in particular, we show that their average intersection numbers behave in a different way when $N$ tends to infinity. For the $W_2^1$ norm, the average intersection numbers of curves of degree $N$ grow only as  $ N$, and for $W_2^r$ with $r>1$ they tend to a finite limit.   

\begin{remark} \rm
Another natural problem is to calculate the average numbers of pairs, triples, ... $k$-ples... of intersection points of pairs of curves corresponding to points of $B^2(N,R)$. These numbers are the coefficients of the Taylor expansion at $t=1$ of the generating polynomial $A_0(N)+A_1(N) t + A_2(N) t^2 + \dots$, in particular, they allow one to reconstruct all the numbers $A_i(N)$. 
Probably the exact calculation of these numbers is difficult, but the asymptotics at $N \to \infty$ might be easier; I also hope that our present approach (with integration over configuration spaces) can simplify these calculations.

This work is motivated by the problem of calculating (and adequately defining) the relative measures of the sets of links with fixed topological types in $\RR^3$.  
\end{remark}

\section{The first reduction}
\label{fred}

Let $T^2$ be the product of two source circles of our maps $(f, g): S^1 \sqcup S^1 \to \RR^2$. Let $\varphi, \psi \in \RR/(2\pi \ZZ)$
be the standard coordinates in these circles, so that $f(\varphi)=(x_1(\varphi),y_1(\varphi))$ and $g(\psi)=(x_2(\psi),y_2(\psi))$.
Consider the space $B^2(N,1) \times T^2$ and its submanifold $I(N)$ consisting of pairs 
\begin{equation}
\label{mainpair}
\left( (f,g),(\varphi,\psi)\right)
\end{equation} 
such that $f(\varphi)=g(\psi)$. It is the space of a fiber bundle over $T^2$, whose fiber over a point $(\varphi,\psi)$ is the subset $L(\varphi,\psi) \subset B^2(N,1)$ consisting of the pairs $(f,g)$ such that $f(\varphi)=g(\psi)$. This subset is a submanifold of codimension 2 in the Euclidean space $\RR^{4(2N+1)}$, in particular, it has the standard volume form. We also have the standard volume form $d\varphi d\psi$ on the torus $T^2$, so that the volume of this torus is $(2\pi)^2$. Having the volume forms on the base and on the fibers of the fiber bundle $I(N) \to T^2$, we also obtain their product volume form on this manifold $I(N)$.

\begin{lemma}
\label{forg}
The differential of the forgetful map $I(N) \to B^2(N,1)$ $($i.e. of the map sending any pair $((f,g)(\varphi_0,\psi_0))$ to the point $(f,g))$ at
any point $((f,g)(\varphi_0,\psi_0))$ 
multiplies the volumes by     
\begin{equation}
\label{forgg}
\frac{1}{2N+1} \left| \det \left|\frac{d f}{d \varphi}(\varphi_0), \frac{d g}{d \psi}(\psi_0) \right| \right|.
\end{equation}
\end{lemma} 

{\it Proof.} This statement is invariant under rotations of the source circles, therefore it
is enough to prove it at any point $((f,g),(\varphi_0,\psi_0))$ with $\varphi_0=0, \psi_0=0$.
Let our two curves $f, g: S^1 \to \RR^2$ be given by (\ref{on}) and (\ref{twoo}).
Then the plane $L(0,0)$ is defined by the system of equations 
\begin{equation}
\label{coo}
\left\{
\begin{array}{ccl}
a_0+ \dots + a_N & = & c_0 + \dots  + c_N , \\  
\tilde a_0+ \dots + \tilde a_N & = & \tilde c_0 + \dots +\tilde c_N  \ .
\end{array}
\right. 
\end{equation}
This plane is orthogonal to the two mutually orthogonal vectors in $B^2(N,\infty)$, the first of which is given by 
\begin{equation}
\label{ort1}
\left\{\begin{array}{cclc}
x_1(\varphi) & = & \frac{1}{2}+\cos \varphi +\dots+\cos N\varphi, & y_1(\varphi)=0, \\
 x_2(\psi) & =& -\left(\frac{1}{2}+\cos \psi + \dots + \cos N\psi \right), & y_2(\psi)=0, 
\end{array}
\right. 
\end{equation}
and the second one by
 \begin{equation}
\label{ort2}
\left\{
\begin{array}{cccl}
x_1(\varphi)=0, & y_1(\varphi) & =& \frac{1}{2}+\cos \varphi +\dots+\cos N\varphi, \\ 
x_2(\psi)=0, &  y_2(\psi) & = & - \left(\frac{1}{2}+\cos \psi + \dots + \cos N\psi \right) \ .
\end{array}
\right. 
\end{equation}

Denote by $X$ and $Y$ two tangent vectors at the point $((f,g),(0,0)) \in B^2(N,1) \times T^2$ whose projections to $B^2(N,1)$ coincide with these two vectors (\ref{ort1}), (\ref{ort2}), and projections to the tangent space of $T^2$ are equal to $0$. Then the vectors 
\begin{equation}
\label{tang1}
\frac{\partial}{\partial \varphi} - \frac{1}{2N+1}\left(\frac{d x_1}{d \varphi}(0) X + \frac{d y_1}{d \varphi}(0) Y\right)
\end{equation} 
and  
\begin{equation}
\label{tang2}
\frac{\partial}{\partial \psi} + \frac{1}{2N+1}\left(\frac{d x_2}{d \psi}(0) X + \frac{d y_2}{d \psi}(0) Y\right)
\end{equation}
are orthogonal to the subspace $L(0,0) \times (0,0)$ in $B^2(N,1) \times T^2$ and are tangent to the manifold $I(N)$: indeed, the family of pairs of maps 
$$(f_\tau,g_\tau) \equiv (f,g) - \tau  \frac{1}{2N+1}\left(\frac{d x_1}{d \varphi}(0) X + \frac{d y_1}{d \varphi}(0) Y\right)$$
parameterized by $\tau \in (-\pi,\pi)$
satisfies the estimate $|f_\tau(\tau) - g_\tau(0)| = O_{\tau \to 0}(|\tau|^2),$ and the family of pairs of maps $$(f_\varkappa,g_\varkappa) \equiv (f,g) + \varkappa  \frac{1}{2N+1}\left(\frac{d x_2}{d \psi}(0) X + \frac{d y_2}{d \psi}(0) Y\right)$$
satisfies $|f_\varkappa(0) - g_\varkappa(\varkappa)| = O_{\varkappa \to 0}(|\varkappa|^2).$

Now let $\Box$ be a parallelepiped of Euclidean volume 1 in $L(0,0)$ (more precisely, in the tangent space of $L(0,0)$). 
The parallelepiped in the tangent space of $I(N)$ at the point $((f,g),(0,0))$ spanned by $\Box$ and two vectors (\ref{tang1}), (\ref{tang2}), has volume 1 in the above-described measure on $I(N)$. The differential of our forgetful map sends this parallelepiped to the parallelepiped in (the tangent space of) $B^2(N,1)$, which is the product of $\Box$ and a parallelogram orthogonal to $L(0,0)$ and spanned by the vectors  $$ \frac{-1}{2N+1}\left(\frac{d x_1}{d \varphi}(0) X + \frac{d y_1}{d \varphi}(0) Y\right) \ \ \mbox{and} \ \  \frac{1}{2N+1}\left(\frac{d x_2}{d \psi}(0) X + \frac{d y_2}{d \psi}(0) Y\right).$$ The Euclidean area of the latter parallelogram is  

$$\frac{1}{(2N+1)^2}  \left| \left|
\begin{array}{cc}
\frac{d x_1}{d \varphi}(0) & \frac{d x_2}{d \psi}(0) \\
\frac{d y_1}{d \varphi}(0) & \frac{d y_2}{d \psi}(0)
\end{array}
\right| \right| \|X\| \|Y\| =  \frac{1}{2N+1}\left| \det \left|\frac{d f}{d \varphi}(0), \frac{d g}{d \psi}(0) \right| \right| . \Box    $$

\begin{corollary}
\label{maincorr} Number 
$$\sum_{i=0}^\infty i \cdot \tilde A_i(N,1)$$ is equal to the integral of the function $($\ref{forgg}$)$ along the manifold $I(N)$.
\hfill $\Box$
\end{corollary}

On the other hand, this integral is equal to $4 \pi^2$ times the integral of the same function along an arbitrary fiber $L(\varphi, \psi)$
with respect to its Euclidean measure. Indeed, it is obvious that all last integrals corresponding to different $\varphi$ and $\psi$ are equal to one another. In the next three sections, we will calculate such an integral for the plane $L(0,0)$, which will give us a proof of Theorem \ref{mth}. \medskip

\section{Integral along the disc}

\begin{lemma}
\label{lem2a}
The images of the points of $S^1$ under the maps $f \in B(N,1)$ can be any points of $\RR^2$ with the norm $\sqrt{x^2+y^2} \leq \sqrt{(2N+1)/2}.$
\end{lemma}

{\it Proof.} The set of such images is obviously invariant under rotations of $\RR^2$. The maximal value of the $x$ coordinate of such an image is equal to the maximal possible value of the function $a_0+ a_1 + \cdots + a_N$ on the domain in $\RR^{N+1}$ given by the inequality $2a_0^2+ \sum_{j=1}^N a_j^2 \leq 1.$  \hfill $\Box$ \medskip

So let's denote by $D(N)$ the disc of radius $\sqrt{(2N+1)/2}$ in $\RR^2$.

Consider the linear map \begin{equation}
\label{secbund}
L(0,0) \to D(N)
\end{equation}
 assigning to any pair $(f,g) \in L(0,0)$ the value $f(0)\equiv g(0)$.
There are two volume measures on $L(0,0)$: one is the Euclidean measure induced from $B^2(N,1)$, and the other is the product of the similar Euclidean measure along the fibers of the map (\ref{secbund}) and the  measure lifted from the standard volume form $dx dy$ on the base of this bundle.

\begin{lemma}
\label{proj}
The former measure is equal to $\frac{4}{2N+1}$ times the latter one.
\end{lemma} 

{\it Proof.} All fibers of the bundle (\ref{secbund}) are orthogonal to two mutually orthogonal vectors in $B^2(N,\infty),$ the first of which is given by 
$$\left\{
\begin{array}{cclc}
x_1(\varphi) & =& \frac{1}{2}+\cos \varphi +\dots+\cos N\varphi, & y_1(\varphi)=0, \\
x_2(\psi) & =& \frac{1}{2}+\cos \psi + \dots + \cos N\psi, & y_2(\psi)=0 \ ,
\end{array}
\right.
$$ 
and the second one by
$$
\left\{
\begin{array}{cccl}
x_1(\varphi)=0, & y_1(\varphi) & = & \frac{1}{2}+\cos \varphi +\dots+\cos N\varphi, \\ 
x_2(\psi)=0, & y_2(\psi) & = & \frac{1}{2}+\cos \psi + \dots + \cos N\psi \ .
\end{array}
\right.
$$ 
The norms of these vectors are equal to $\sqrt{2N+1}$, 
and their projections to $D(N)$ are 
vectors $(\frac{1}{2}+N) \partial/\partial x$ and $(\frac{1}{2}+N) \partial/\partial y.$ Therefore, the ratio  of the areas of squares spanned by these pairs of vectors is equal to $$\hspace{3.8cm} \frac{\sqrt{2N+1}^2}{(\frac{1}{2}+N)^2} \equiv \frac{4}{2N+1} \ .  \hspace{3.8cm} \Box $$ 

\medskip

Given a point $U \in D(N)$, consider the subset $\Lambda(U) \subset B(N,1)$ consisting of all maps $f:S^1 \to \RR^2,$ $f \in B(N,1)$, such that $f(0)=U$. It is the intersection of the ball $B(N,1)$ and an affine plane of codimension 2 in the Euclidean space $\RR^{4N+2}$ of all maps (\ref{on}), therefore also a ball. 
The fiber of the bundle (\ref{secbund}) over the point $U$  is the product of this ball $\Lambda(U)$ and the analogous ball in the space of maps $g$. 
Define the number $\tilde \Xi(U)$ as the integral  of the function $|f'(0)|$ along this ball $\Lambda(U)$ (with respect to its Euclidean measure). This number obviously depends only on the norm of $U$, so we can define $\Xi(A)$, $A \in \left[0,\sqrt{(2N+1)/2}\right],$ as the common value of $\tilde \Xi(U)$ for all $U$ with $|U|=A$. 

\begin{prop}
\label{prop3}
1. The integral of the function $\left| \det \left|\frac{d f}{d \varphi}(0), \frac{d g}{d \psi}(0) \right| \right| $ along the squared ball $\Lambda(U) \times \Lambda(U) \subset L(0,0)$ $($that is, along the fiber of the map $($\ref{secbund}$)$ over the point $U)$ is equal to 
$
\frac{2}{\pi} (\tilde \Xi(U))^2.
$

2. The integral of the function $($\ref{forgg}$)$ over the manifold $L(0,0)$ is equal to
\begin{equation}
\label{need}
\frac{1}{2N+1} \cdot \frac{2}{\pi} \cdot \frac{4}{2N+1} \int_0^{\sqrt{(2N+1)/2}} 2\pi A \cdot (\Xi(A))^2 dA.
\end{equation}
\end{prop}

{\it Proof.} 1. It suffices to prove this for $U$ lying in the $x$-axis, $U= (A,0)$. In this case
 circle $S^1 \sim SO(2,\RR)$ acts on $\Lambda(U)$ by linear isometries, preserving all  coefficients $a_j$ and $\tilde a_j$ in (\ref{on}), and moving coefficients $b_j$ and $\tilde b_j$ by the rotations of the plane $\RR^2$: namely, for any $j=1, \dots, N$, the point $\theta \in [0,2\pi]$ of the acting circle moves the numbers $b_j$ and $\tilde b_j$ to $b_j \cos \theta - \tilde b_j \sin \theta$ and $b_j \sin \theta+\tilde b_j \cos \theta$ respectively. 

Consider also the action of $SO(2,\RR)$ on $\Lambda(U) \times \Lambda(U)$, which is the product of the previous one on the first factor $\Lambda(U)$ and the trivial action on the second factor. The  number $(\tilde \Xi(U))^2$ can be considered as the integral of the function $|f'(0)| \cdot |g'(0)|$ along the domain $\Lambda(U) \times \Lambda(U)$. The integral of the function $\left| \det \left|\frac{d f}{d \varphi}(0), \frac{d g}{d \psi}(0) \right| \right| \equiv |f'(0)|\cdot |g'(0)| \cdot |\sin\angle(f'(0),g'(0))| $ over any orbit of our $SO(2,\RR)$-action is equal to $\frac{2}{\pi}$ times the integral of the function $|f'(0)| \cdot |g'(0)|$ over the same orbit; here the factor $\frac{2}{\pi}$ is the Buffon constant, that is, the average value of $|\sin \theta|$ on $[0,2\pi]$, see \cite{San}.

2. Statement 2 immediately follows from statement 1 and Lemma \ref{proj} (providing the factor $\frac{4}{2N+1}$ in (\ref{need})) by integration over the disc $D(N)$ in polar coordinates. \hfill $\Box$

\section{Calculation of the integral $\Xi(A)$}

 We again assume  that $U$ belongs to the $x$-axis, $U=(A,0)$, $A \geq 0$. 
Then our $4N$-dimensional ball $\Lambda(U) \subset B(N,1)$ is distinguished by the equations 
\begin{equation}
\label{dis}
\left\{
\begin{array}{ccc}
a_0+ a_1+ \cdots + a_N & = &A, \\
\tilde a_0 + \tilde a_1 + \cdots + \tilde a_N & = & 0.
\end{array}
\right. 
\end{equation}
without additional restrictions on the coordinates $b_j, \tilde b_j$. 
Its closest point to the origin has coordinates 
$$a_0=\frac{A}{2N+1}, \  a_1= \cdots = a_N=\frac{2A}{2N+1}, $$ with all $\tilde a_j, b_j$ and $\tilde b_j$ equal to 0. The norm of this point is equal to $\sqrt{2/(2N+1)} A.$ Accordingly, the radius of the ball $\Lambda(U)$ is equal to $\sqrt{1-\frac{2A^2}{2N+1}}$. Denote this number by $k(A)$.
We have two linear functions on this ball, $\beta \equiv \sum_{j=1}^N j b_j$ and $\tilde \beta = \sum_{j=1}^N j\tilde b_j$, 
which are equal to the derivatives $\frac{d x}{d\varphi}(0)$ and $\frac{dy}{d\varphi}(0)$ of the map (\ref{on}). The norms of these linear functions in the dual Euclidean space are equal to $\sqrt{1+4+ \dots + N^2}$; let's denote this number by  $\lambda(N)$.
By an orthogonal change of coordinates, we can turn these functions to $\lambda(N) b$ and $\lambda(N) \tilde b$ respectively, where $b$ and $ \tilde b$ are some new Euclidean coordinates. 
Consider the orthogonal projection of our ball $\Lambda(U)$ to the coordinate plane $\RR^2$ of these two coordinates. Its image is the disc of radius $k(A)$, and the pre-image of any point with norm $r$ is a $(4N-2)$-dimensional ball of radius $\sqrt{k(A)^2 -r^2}.$ Function \begin{equation}
\label{veloc}
|f'(0)| \equiv \sqrt{\left(\sum_{j=1}^N j b_j\right)^2 + \left(\sum_{j=1}^N j \tilde b_j\right)^2}
\end{equation}
 is constant on these fibers, namely, it is equal to $\lambda(N)r$.

Therefore, our integral $\Xi(A)$ is equal to
$$\int_{0}^{k(A)}  2\pi r \cdot \lambda(N) r \cdot \frac{\pi^{2N-1}}{(2N-1)!} \sqrt{k(A)^2-r^2}^{4N-2} dr,$$
which by substitution \
$r=k(A) \sin \alpha$ \ equals
$$\frac{2\pi^{2N}}{(2N-1)!} (k(A))^{4N+1} \lambda(N) \int_0^{\pi/2} \sin^2 \alpha \cos^{4N-1} \alpha \ d\alpha=$$
\begin{equation}
\label{xiansw}
=\frac{\pi^{2N}2^{4N}(k(A))^{4N+1} \lambda(N) \cdot (2N)!}{(4N+1)!} \ .
\end{equation}

\section{Calculation of integral (\ref{need})} 
\label{calin}

Substituting (\ref{xiansw}) to (\ref{need}) we get 
$$\frac{16}{(2N+1)^2}\frac{\pi^{4N} 4^{4N}(\lambda(N))^2 ((2N)!)^2}{((4N+1)!)^2} 
\int_0^{\sqrt{(2N+1)/2}} A \cdot (k(A))^{8N+2} dA .$$

The substitution $A=\sqrt{(2N+1)/2} \sin \gamma$ turns the last integral to
\begin{equation}
\label{eight} 
\frac{2N+1}{2}  \int_0^{\pi/2}  \sin \gamma (\cos \gamma)^{8N+3} d \gamma = \frac{1}{8},
\end{equation}
and so the number (\ref{need}) is 
$$  \frac{2}{(2N+1)^2}\frac{\pi^{4N} 4^{4N}(\lambda(N))^2 ((2N)!)^2}{((4N+1)!)^2} \ .    $$
It remains to multiply it by $4\pi^2$ (the volume of the configuration space of points $(\varphi,\psi) \in T^2$) and divide by the volume (\ref{totvol}) of the domain $B^2(N,1)$ to get the formula
(\ref{mfor}). \hfill $\Box$

\section{The same for $W_2^r$-metric}
\label{sob}

Let  now the space of curves (\ref{on}) be equipped with  
$W_2^r$ Euclidean structure: the dot product of two curves $(x_1,y_1)$ and $(x_2,y_2)$ equals 
$$\frac{1}{\pi}\int_{0}^{2\pi} \sum_{q=0}^r \left(x_1^{(q)}(\varphi)x_2^{(q)}(\varphi) + y_1^{(q)}(\varphi)y_2^{(q)}(\varphi)\right ) d \varphi $$ 
where $^{(q)}$ means the $q$-th derivative.
In particular, the scalar product of curves (\ref{on}) and (\ref{twoo}) is equal to 
$$2\left(a_0c_0 +\tilde a_0 \tilde c_0\right)+ \sum_{j=1}^N\left(\left(1+j^2+ \cdots + j^{2r}\right)\left(a_j c_j+ b_j d_j+\tilde a_j \tilde c_j + \tilde b_j \tilde d_j\right)\right) .$$
 
The definition of this space and its Euclidean structure is also invariant under the rotations of the plane $\RR^2$.

\begin{definition} \rm
Denote the number $1+j^2+ \cdots + j^{2r}$ by $\tau_j$, and the number $1+2\sum_{j=1}^N \frac{1}{\tau_j}$ by $\mu(N)$; define also the number 
\begin{equation}
\label{lamr}
\lambda_r(N) \equiv \left(\sum_{j=1}^N \frac{j^2}{\tau_j} \right)^{1/2} 
\end{equation}
 (that is, the maximal value of the function $|\sum_{j=1}^N j \cdot b_j|$ on the ellipsoid $\left\{ \sum_{j=1}^N (\tau_j \cdot b^2_j) \leq 1 \right\} $ ).
\end{definition}

\begin{theorem}
\label{mth2}
For any non-negative integer $r$, the average number of intersection points of two curves of the form $($\ref{on}$)$ from the unit ball in the  $W_2^r$ metric is equal to
\begin{equation}
\label{eqmn2}
\frac{2^{8N+3} (\lambda_r(N))^2  ((2N)!)^4 (2N+1)} {\mu(N) ((4N+1)!)^2} \ . 
\end{equation}
\end{theorem}

\begin{corollary}
If $r=1$ then this number behaves asymptotically as $ c\cdot N(1+ o_{N \to \infty}(1)) $ for some positive constant $c$.
For any natural $r\geq 2$ this number tends to a finite non-zero limit as $N \to \infty.$ These limits decrease as $\frac{2\pi}{r^2 }(1+ o(1))$ when $r$ tends to infinity. \hfill $\Box$
\end{corollary}

{\it Proof} of Theorem \ref{mth2} is very similar to that of Theorem \ref{mth}; throughout this proof the number $\mu(N)$ replaces $2N+1$. In particular, in the analog of Lemma \ref{forg}, the factor $\frac{1}{2N+1}$ in (\ref{forgg}) should be replaced by $\frac{1}{\mu(N)}$, and the factor $\frac{4}{2N+1}$ in Lemma \ref{proj} should be replaced by $\frac{4}{\mu(N)}$. 
The analog of the disc $D(N)$ of values of the curves from the unit ball $B(N,1)$ has now the radius equal to $\sqrt{\frac{\mu(N)}{2}} ,$
and the radius $k(A)$ of the ball $\Lambda(U)$, $|U|=A$, is equal to $\sqrt{1-\frac{2A^2}{\mu(N)}} .$ 

The analog of the number $\lambda(N) \equiv \sqrt{1+4+\cdots+N^2}$ is now called $\lambda_r(N)$, see (\ref{lamr}).
The rest of the proof is completely the same as that of Theorem \ref{mth}. (Attention: the analog of the expression in the left-hand part of (\ref{eight})
in this calculation is equal not to $1/8$ but to $\frac{\mu(N)}{8(2N+1)}$ : its simplification in the case $r=0$ is an occasional coincidence.)  \hfill $\Box$

\end{document}